\documentclass{amsart}
\usepackage{amsmath}
\usepackage{amsthm}
\usepackage{amssymb}
\usepackage{tikz-cd}
\usepackage{courier}
\usepackage{amsfonts}
\usepackage{multicol}
\usepackage[mathscr]{euscript}
\usepackage{dirtytalk}
\usepackage{graphicx}
\graphicspath{{./images/}}
\newcommand{\interior}[1]{%
  {\kern0pt#1}^{\mathrm{o}}%
  }
\newtheorem{theorem}{Theorem}[section]
\newtheorem{lemma}[theorem]{Lemma}
\newtheorem{corollary}[theorem]{Corollary}

\newtheorem{definition}{Definition}[section]
\newtheorem{conjecture}[theorem]{Conjecture}

\newtheorem*{theorem*}{Theorem}
\newtheorem{remark}[theorem]{Remark}
\begin{document}
\title{Geometry of knots in real projective 3-space}
\author{Rama Mishra}
\address{Department of Mathematics, Indian Institute of Science Education and Research, Pune, India}
\email{r.mishra64@iiserpune.ac.in}

\author{Visakh Narayanan}
\address{Department of Mathematics, Indian Institute of Science Education and Research, Pune, India}
\email{vishakhme@gmail.com}
\maketitle
\begin{abstract}
This paper discusses some geometric ideas associated with  knots in  real projective 3-space $\mathbb{R}P^3$. These ideas are borrowed from classical knot theory.  Since knots in $\mathbb{R}P^3$ are classified into three disjoint classes, - affine, class-$0$ non-affine and class-$1$ knots, it is natural to wonder in which class a given knot belongs to. In this paper we attempt to answer this question. We provide
 a structure theorem for these knots  which helps in describing their behaviour near the projective plane at infinity. 
 We propose a procedure called {\it space bending surgery}, on affine knots to produce several examples of knots.  We later show that this operation  can be extended on an arbitrary knot in $\mathbb{R}P^3$. We also define a notion of \say{ genus} for knots in $\mathbb{R}P^3$ and study some of its properties. We prove that this genus detects knottedness in $\mathbb{R}P^3$ and gives some criteria for a knot to be affine and of class-$1$. We also prove a \say{non-cancellation} theorem for space bending surgery using the properties of genus. We produce examples of class-$0 $ non-affine knots with genus $1$. And finally we study the notion of companionship of knots in $\mathbb{R}P^3$ and using that we provide a geometric criteria for a knot to be affine. 
 Thus we highlight that, $\mathbb{R}P^3$ admits a knot theory with a truly different flavour than that of $S^3$ or $\mathbb{R}^3$. 
\end{abstract}

\section{Introduction}

A knot, in a smooth manifold $M$ is an embedding $K:S^1\to M$ which is recogonized upto ambient isotopy of $M$ (\cite{crow}). As usual, most of the times, we would not distinguish between the embedding and its image. 
Classically, knot theory was studied in $S^3$ or $\mathbb{R}^3$. Here we study the knot theory of $\mathbb{R}P^3$.
 The canonical two sheeted covering map,
$h: S^3\to \mathbb{R}P^3$
is an essential tool for understanding the geometry of $\mathbb{R}P^3$. We will denote this map by $h$ throughout in this manuscript.
The antipode map $a:S^3\to S^3$, is the only non-trivial automorphism of the covering map. The fundamental group of $\mathbb{R}P^3$ is $\frac{\mathbb{Z}}{2\mathbb{Z}}$.  Thus, a  knot in $\mathbb{R}P^3$ can represent either $\bar{0}$ or $\bar{1}$ in the fundamental group. We refer to knots in $\mathbb{R}P^3$ to be of class-$0$ or of class-$1$ accordingly. All class-$0$ knots, will lift to a link with two components in $S^3$. The two components are antipodes of each other. The class-$1$ knots lift to two paths which together form a knot in $S^3$. This knot is the antipode of itself. 
Of late  there has been a lot of interest in studying knots in $\mathbb{R}P^3$ and their lifts in $S^3$ (\cite{enrico}). Also some of the knot invariants from classical knot theory  have been generalized. For example, in \cite{julia}  an analogue of the Jones polynomial for links in $\mathbb{R}P^3$ and in \cite{ Huy-Le} Twisted Alexander polynomial have been developed. In this paper, we focus  more on the geometric properties associated with these knots.

\vspace{.1in}

Many of the tools used to study classical knots are not available in $\mathbb{R}P^3$. For instance, one of the techniques for constructing knots  in $S^3$ is the notion of connected sum for knots that gives rise to a notion of primeness of a knot. Recall that the connected sum of knots $K_1$ and $K_2$ in $S^3$  is defined using a surgery procedure as :    Think of $K_1$ and $K_2$ to be contained in two distinct copies of $S^3$. Now remove an open ball $B_1$ from the sphere containing $K_1$, which intersects $K_1$ at a single arc with its two boundary points on $\partial \bar{B_1}$. Similarly remove an open ball $B_2$ from the sphere containing $K_2$, which intersect $K_2$ at a single arc. The complements of $B_1$ and $B_2$ are closed 3-balls. By gluing their boundaries together, with a diffeomorphism which also sends the boundary points of the remaining arc of $K_1$ to the boundary points of the remaining arc of  $K_2$, we obtain a copy of $S^3$ containing a new knot.  Note that the ambiguity concerning the choice of points in $\partial (K_1\setminus B_1)$ to be identified with specific points in $\partial( K_2\setminus B_2)$  can be taken care   if we start with an orientation on both the knots $K_1$ and $K_2$ \cite{rolfsen}. This gives a well defined notion of $K_1\#K_2$ as a new knot in $S^3$. This procedure can be performed repeatedly as a sum of an arbitrary number of knots and it is an associative binary operation which defines a monoid structure on the set of knots in $S^3$. Unfortunately this process is not well defined in $\mathbb{R}P^3$. 
  
\par In this paper, we develop a surgery procedure which is a method of generating new knots in $\mathbb{R}P^3$ from a given knot. In  some sense it is described as a notion of adding a projective line to the knot. This procedure can be repeated arbitrarily many times and each time we perform it, the homology class of the knot strictly changes. We refer to it as \say{class changing surgery}.

\par Another property which projective knots don't have and classical knots have is the existence of a Seifert surface. It is easy to see that if a knot in any three manifold $M$, is the boundary of an orientable surface, then the knot has to be homologically trivial in $M$. Hence for any three manifold which has a non-trivial first homology, the knots in it cannot have a complete theory of Seifert surfaces. For example the class-1 knots in $\mathbb{R}P^3$ cannot have Seifert surfaces. Thus one cannot directly generalize the notion of genus for projective knots. In this paper we define a notion of \say{genus} for a knot in $\mathbb{R}P^3$. The philosophy we adopt here, is that,  \say{the complexity of a knot should also reflect in the complexity of the surfaces on which it can be embedded}. The properties of this genus is very different from that of the classical genus. Although there are some similarities too. For instance this genus detects knottedness similar to the genus. It also gives a criteria for a knot to be affine or to be of class-1. We also construct several examples of class-0 non-affine knots of genus 1. But we do not completely understand how the genus of this class behaves. And thereby it reveals that the class-0 non-affine knots is a tough family of knots.  
\par We also generalize the notion of companionship of knots from classical knot theory to $\mathbb{R}P^3$. It is known that the classical unknot is a companion of every knot in $S^3$. As a parallel to this, we prove that the projective unknot is a companion to every knot in $\mathbb{R}P^3$. And finally we  provide a geometric characterization for affineness of a knot which we call as \say{projective inessentiality}. 

\textbf{Organization of the paper:} In Section 2, we define the ball model of $\mathbb{R}P^3$ and the notion of diagram of a knot in $\mathbb{R}P^3$. These ideas arise frequently in the paper, hence most of the notations used in this section has been used throughout the paper. Then in Secion 3, the terms residual tangles, space bending surgery and class changing surgery have been introduced using which a structure theorem for projective knots has been proved. In Section 4, the genus for projective knots is defined and some of its properties are discussed. Section 5 defines the  concept of companionship for knots in $\mathbb{R}P^3$ and discusses its comparison with the existing notion of companionship for classical knots. At the end of Section 5 the notion of projective inessentiality is defined and it has been shown that it provides a complete criterion for a given knot in $\mathbb{R}P^3$ to be affine.

\section{The geometry and topology of $\mathbb{R}P^3$}
Topology of $\mathbb{R}P^3$ can be understood in many ways. One way is to think of it as  a compactification of $\mathbb{R}^3$ by adding a projective plane at infinity. Adding an $\mathbb{R}P^2$ at infinity is done in the following way. The complement of an open $3$-ball in $\mathbb{R}^3$ is homeomorphic to $S^2\times [0,1)$. Choose a closed 3-ball $B\subset \mathbb{R}^3$. Choose a tubular neighborhood of $\partial B$ and remove it from $\mathbb{R}^3$. We may identify the remaining two components with $B$ and $S^2\times [0,1)$. First we compactify this $S^2\times [0,1)$ to $S^2\times [0,1]$ by adding a $2$-sphere (at infinity) as $S^2\times\{1\}$. Then we form a quotient of the $S^2\times[0,1]$ by identifying the points of the form $(p,1)$ to the points of the form $(-p,1)$. Let $C$ be space thus obtained. Notice that the quotient map takes the sphere $S^2\times\{1\}$ in $S^2\times [0,1]$ and projects it onto a projective plane in $C$. The boundary of $C$ is a $S^2$ which is the image of $S^2\times\{0\}$ under the quotient map. Now by gluing the boundaries of $B$ and $C$, we obtain a copy of $\mathbb{R}P^3$. Notice the complement of any $3$-ball in $\mathbb{R}P^3$ is homeomorphic to $C$. Also by construction, it must be clear that $C$ is homeomorphic to the mapping cylinder \cite{hatcher} of the two sheeted covering $S^2\to \mathbb{R}P^2$. Thus $\mathbb{R}P^3$ can also be thought of as,  joining a $3$-ball  to this mapping cylinder along their boundary sphere. We will use this decomposition frequently.\\ 
\par  $\mathbb{R}^3$ is homeomorphic to an open ball. The process we described above, can alternatively be described as follows. Start with a closed $3$-ball and identify the diagonally opposite points in its boundary sphere. Thus we obtain $\mathbb{R}P^3$ as a quotient of the 3-ball. 
\begin{definition}
A description of $\mathbb{R}P^3$ as the quotient of a 3-ball will be called as a \say{ball model} for $\mathbb{R}P^3$.   
\end{definition}
 All lines and planes in $\mathbb{R}^3$ get compactified as projective lines and projective planes respectively in $\mathbb{R}P^3$. The complement of any projective plane, which is a compactification of a plane in $\mathbb{R}^3$, in $\mathbb{R}P^3$ is diffeomorphic to $\mathbb{R}^3$. The rich geometry admitted by $\mathbb{R}P^3$, similar to the euclidean space, can be exploited to undertand the knots in it. For example, the theory of diagrams for knots in $S^3$ as in, for more details see \cite{crow}, can be generalized to knots in $\mathbb{R}P^3$ as in \cite{viro}. Since this notion of diagrams will be used throughout the manuscript, we wish to give a quick exposition of this here. Most of the notations used here for certain functions or spaces will be used later without defining them again. Since these notations are not used here for any other purpose, there should not be any confusion arising from this.   
\par Corresponding to any line $L$ in $\mathbb{R}^3$, we can construct a family of lines consisting only of all parallel lines to $L$. Every point in $\mathbb{R}^3$ belongs to exactly one line in this family. Now corresponding to any such family there exists a plane $Z$ passing through origin, such that, $Z$ intersects each line parallel to $L$ at a unique point. Thus, we can define a projection, $d:\mathbb{R}^3\to Z$, sending a point of $\mathbb{R}^3$ lying on a line $L'$ parallel to $L$ to the unique point in $L'\cap Z$. Notice that there is a unique line in the family which passes through origin. We may regard $L$ itself as this line, without loss of generality. Now $L$ and $Z$ are also vector subspaces of the vector space $\mathbb{R}^3$. Thus we get a splitting $\mathbb{R}^3\approx L\oplus Z$. It is easy to see that the map $d$ defined above is same as the projection of the vector space $\mathbb{R}^3$ onto its direct summand $Z$ in the previous splitting. By using this projection map one can have a diagram for a knot $K$ in $\mathbb{R}^3$ on the plane $Z$. Note that by using different choices of $L$ and $Z$ we have several such projections.

\par Notice that any family of parallel lines in $\mathbb{R}^3$ meet at a unique point in $\mathbb{R}P^3$ at infinity. Different points in the projective plane at infinity correspond to different family of parallel lines. Consider the ball model of $\mathbb{R}P^3$ presented as a quotient of the unit $3$-ball $D^3\subset \mathbb{R}^3$. Choose two diametrically opposite points $N$ and $N'$. Choose a great circle $\gamma$ on $\partial D^3$ disjoint from $N$ and $N'$. Let $D^2\subset D^3$, be a 2-disk, whose boundary is $\gamma$. Notice that since $N$ and $N'$ are diametrically opposite, they lie on distinct hemispheres in $\partial D^3\setminus \gamma$. Consider a point $p\in D^3\setminus\{N,N'\}$. It follows from elementary euclidean geometry that, either there exist a unique circle in $\mathbb{R}^3$ passing through $N, N'$ and $p$ or they are collinear points. If $N$, $N'$ and $p$ are collinear, let $L_p$ be the unique diametrical line segment from $N$ to $N'$. It passes through $p$. Suppose $N$, $N'$ and $p$ are not collinear. Then as said above, there is a circle, $\alpha_p$, passing through $N, N'$ and $p$. It easy to see that $\alpha_p\setminus\{N,N'\}$ is just two open arcs out of which exactly one contains $p$. Then let $L_p$ be the arc of $\alpha_p$ which contains $p$ with $N$ and $N'$ added as boundary points. Every circle passing through $N$ and $N'$ has to meet $D^2$ atleast at one point. But notice that for all $p\in D^3\setminus\{N,N'\}$ we know that $L_p$ intersects $D^2$ at a unique point. Also if $q$ lies in $L_p$ then $L_q$ is same as $L_p$. Then we have a projection $\delta: D^3\setminus \{N,N\} \to D^2$, defined by sending every point $p$ to the unique point of $L_p\cap D^2$. It is obvious that it is a continuous projection, with the inverse image of every point $p\in D^2$ as $L_p\setminus\{N,N'\}$. Refer to Figure 1. 
\par Under the quotient map $\pi: D^3\to \mathbb{R}P^3$, $N$ and $N'$ maps to the same point $\overline{N}$ in $\mathbb{R}P^3$. And $D^2$ maps to a projective plane $P$. Each of the arcs $L_p$ maps to a projective line, containing $\overline{N}$. Thus $\delta$ defines a map from $\mathbb{R}P^3\setminus \{\overline{N}\}\to P$. We shall refer to this map also as $\delta$. From the context the reader can deduce which of the projections are being specified. It is easy to see that the points $p\in D^3\setminus \partial D^3$, the open ball, under $\delta$ maps to $D^2\setminus \partial D^2$ the interior of $D^2$. We may identify the open ball to $\mathbb{R}^3$, in such a way that the interior of $D^2$ gets identified to plane $Z$ from the previous example. It is easy to see that there is such an identification where the family of open arcs $L_p\setminus \{N,N'\}$ for all $p$'s in the open ball, maps to a family of parallel lines as before. Hence the projection $\delta$ can be thought of as an extension of the projection, $d:\mathbb{R}^3\to Z$, defined above.  

Any link, which is disjoint from $\overline{N}$ can be projected to $P$. It is obvious that one can always assume that any link is disjoint from $\overline{N}$, since if otherwise, there is an isotopy which can be used to remove it from meeting $\overline{N}$. This image of the projection in $P$ can be \say{drawn} on $D^2$ by taking its inverse image in $D^2$. Clearly it will be a collection of arcs in $D^2$ with boundary points in $\partial D^2$ forming a set closed under the antipode map of $\partial D^2$. As in the classical case, one can always assume the link is in \say{regular position} and the image in $D^2$ has only double points as singularities. 

\begin{definition}
The image of a tame link in $\mathbb{R}P^3$, in regular position, under $\delta$, with the \say{over-under} crossing information provided at each double point, is said to be a \textbf{diagram} of the link.
\end{definition}

Just like the ball, the diametrically opposite points of the disk $D^2$ should be regarded as being identified.
In a diagram of a link in $\mathbb{R}P^3$ on a disk $D^2$, any diameter of the disk is a projective line in the quotient $\mathbb{R}P^2$.
Notice that the inverse image of such a line under $\pi$, in $D^3$ together with $N$ and $N'$ added is a $2$-disk transversally intersecting $D^2$. Refer to Figure 3. And in the quotient $\mathbb{R}P^3$ this disk represents a projective plane. Thus diameters of the disk in a diagram can be pulled back to a projective plane. A projective plane embedded in $\mathbb{R}P^3$ whose inverse image under $\pi$ of some ball model is a flat 2-disk or the sphere $\partial D^3$ in $D^3$ will be called a \say{standard projective plane}.
\begin{definition}
An embedded projective plane in $\mathbb{R}P^3$ which is isotopic to a standard projective plane will be called an asymptotic projective plane. 
\end{definition} 

\begin{remark}
The complement of any asymptotic projective plane is homeomorphic to $\mathbb{R}^3$. 
\end{remark}

\begin{definition}
A knot which has a diagram with no crossings will be called an \textbf{unknot}. Similarly a link which has a diagram with no crossings will be called a \textbf{unlink}. 
\end{definition}

\begin{figure}
\begin{center}
\includegraphics{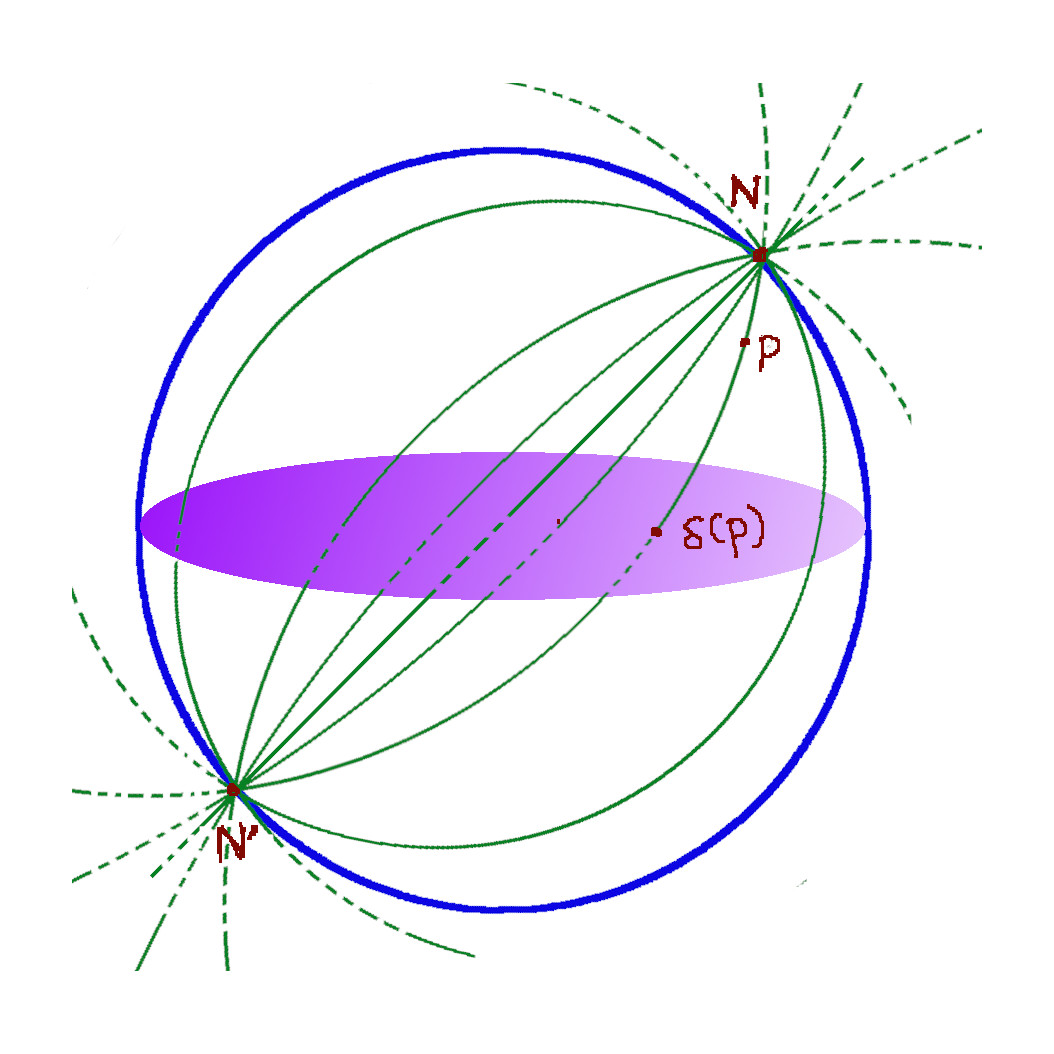}
\end{center}
\caption{The projection map used for making diagrams of knots.}
\end{figure}

 \begin{figure}
\begin{center}
\includegraphics{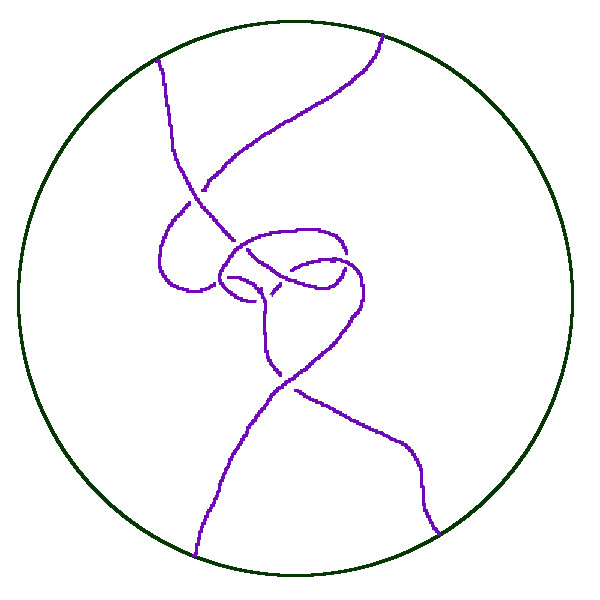}\\
\caption{An example of a knot diagram in $\mathbb{R}P^3$.}
\end{center}
\end{figure}

\begin{figure}
\begin{center}
\includegraphics[scale=1.3]{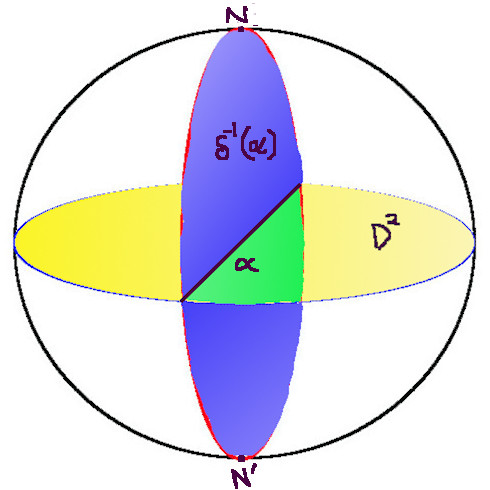}
\caption{Diameters of the target $D^2$ are images of perpendicular disks}
\end{center}
\end{figure}

One of the first differences between classical knot theory and projective knot theory is, there are two different candidates which deserve to be regarded as \say{unknots} in $\mathbb{R}P^3$. We call them  the \say{affine} and \say{projective} unknots. Refer to Figure 4. Notice that any unlink can only have atmost one projective unknot as a component. 
\begin{figure}
\includegraphics{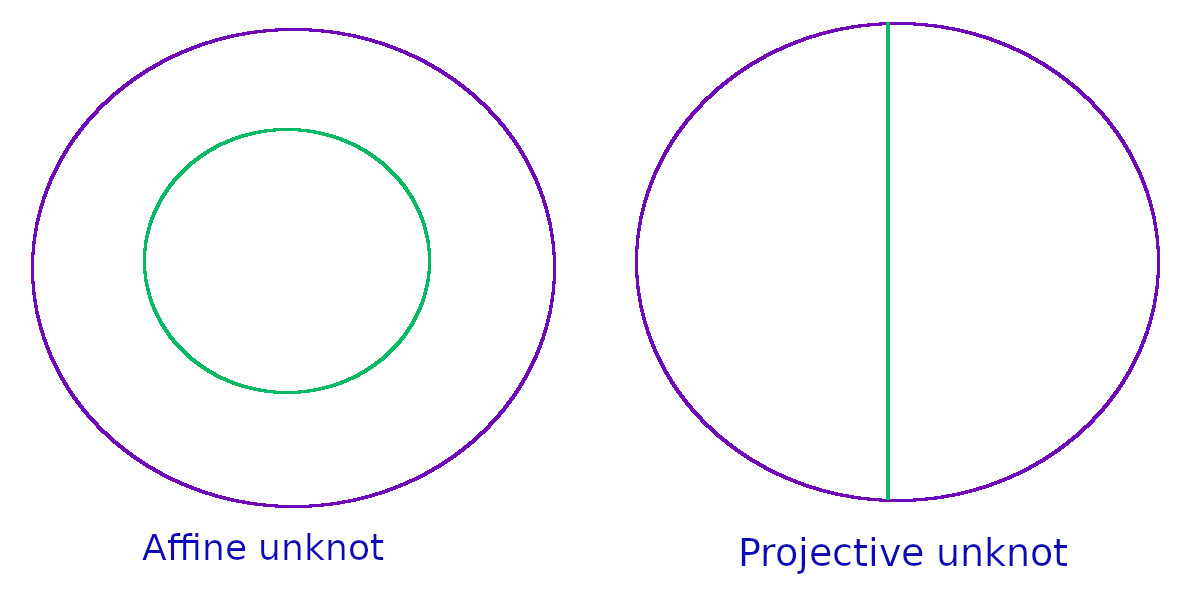}
\caption{Unknots in $\mathbb{R}P^3$}
\end{figure}

\begin{definition}
A knot $K$ is said to be affine, if it is disjoint from some asymptotic projective plane in $\mathbb{R}P^3$. 
\end{definition}

\begin{remark}

\begin{enumerate}

\item Affine knots in $\mathbb{R}P^3$ are thus knots, contained in an open set homeomorphic to $\mathbb{R}^3$. Thus any classical knot, which is a knot in $\mathbb{R}^3$ by definition, has an affine knot in $\mathbb{R}P^3$ representing it. \\
\item Since every knot in $\mathbb{R}^3$ is contractible, the affine knots in $\mathbb{R}P^3$ are of class-$0$. Thus every class-$1 $ knot is non-affine. \\
\item The inverse image of an affine knot is a link in $S^3$ with two components. The inverse image of any asymptotic projective plane in $\mathbb{R}P^3$ is a 2-sphere in $S^3$ whose complement is homeomorphic to two disjoint open 3-balls. Thus for an affine knot $K$, disjoint from an asymptotic projective plane $P$, the two components of its inverse image lie in the two distinct open balls in $S^3$ separated by the 2-sphere lying over $P$. Hence, the linking number of the inverse image of $K$ in $S^3$  must be zero.

\end{enumerate}

\end{remark}

 Now suppose the lift of a class-$0 $ knot $K$ is a link in $S^3$ with a non-zero linking number then clearly $K$ is not an affine knot. The knot shown in Figure 2 is an example where the inverse image has a non-zero linking number. Thus we have three distinct sets of knots to study in $\mathbb{R}P^3$, \textbf{affine, class-0 non-affine and class-1} knots.  

One simple criterion for detecting affineness from a diagram is the following.
\begin{remark}
If $K$ has a diagram on a disk which has a diameter disjoint from the diagram, then $K$ is an affine knot.  
\end{remark}
\begin{figure}
\begin{center}
\includegraphics{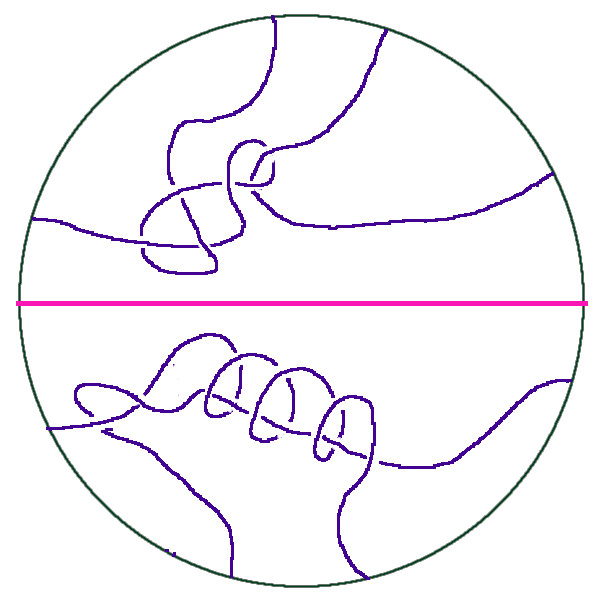}
\caption{Diameter disjoint from the diagram}
\end{center}
\end{figure}
\textbf{Proof of the remark:} Let $d$ be a diameter of the disk $D^2$ on which diagram of $K$ is drawn. Assume, $d$ is disjoint from the diagram. The diameter (as in Figure 5) contains the central point, say $M$, of the 3-disk $D^3$. Notice that $M, N$ and $N'$ are collinear in the ambient $\mathbb{R}^3$ containing $D^3$. Hence there is a plane $Q$ in $\mathbb{R}^3$ containing $N, N'$ and $d$. Then, $Q\cap D^3$ is a 2-disk, $D$ whose quotient in $\mathbb{R}P^3$ is an asymptotic projective plane disjoint from the knot. Then it follows that $K$ is affine.

\section{A structure theorem for knots in $\mathbb{R}P^3$}
\par Consider a smooth 2-sphere $S$ which bounds a ball in $\mathbb{R}P^3$. Removing a tubular neighborhood of $S$ produces two connected components, a $3$-ball $B$ and a mapping cylinder, $C$, of the two sheeted covering map from $S^2$ to $\mathbb{R}P^2$. Since the boundaries of $B$ and $C$ are naturally identified to $S$, for the sake of simplicity, most of the times we will refer to both $\partial B$ and $\partial C$ by $S$. As mentioned above, if we join two knots in $\mathbb{R}P^3$ by removing an unknotted arc contained in a ball and identifying the boundary points as in the classical case, then the resulting knot is in $M:=\mathbb{R}P^3\#\mathbb{R}P^3$.\\
\par The above mentioned process seems to be the only technique to define a connected sum \say{operation} for all knots in $\mathbb{R}P^3$. And it comes with the hassle of dealing with knots in an arbitrary $3$ manifold $M$ and $\mathbb{R}P^3$ at the same time. But for certain special families of knots in $\mathbb{R}P^3$  we can have some operation very similar to the connected sum. Below is a procedure, that may be regarded as a connected sum of an affine knot with the projective unknot.

\subsection{Projectivization of affine knots}

 Consider an affine knot, $K$, in $\mathbb{R}P^3$ which is embedded inside a ball, $B$ with a boundary sphere, $S$. Then we can push a small arc of $K$ out of $B$ in such a way that $K$ now intersects $S$ at exactly two points and the complement $C$ of $B$ intersects $K$ at an unknotted arc disjoint from an asymptotic projective plane in $C$. Then clearly this arc is a trivial cycle in $H_1(C,S)$. Now consider a standardly embedded projective unknot, $J$, in another copy of $\mathbb{R}P^3$. Removing a small ball, $B'$, with a boundary $S'$, containing a unknotted arc of $J$, produces a copy of the mapping cylinder, $C'$. And $J\cap C'$ represents a non-trivial cycle in $H_1(C',S')$. Refer to Figure 8. Now by gluing $B$ and $C'$ along their boundary spheres in such a way that $K\cap S$ gets identified with $J\cap S'$ produces a knot in $\mathbb{R}P^3$. Refer to Figure 6. Notice that the initial knot $K$ is of class-$0$ while the new one is a class-$1$ knot. Thus this is a process of producing a class-$1$ knot from a class-$0$ knot, which we will call \say{projectivization}. A knot in $\mathbb{R}P^3$ which intersects an asymptotic projective plane exactly at one point can be thought of as produced from an affine knot via the above surgery. Hence we wish to make the following definition.
 \begin{figure}
\includegraphics{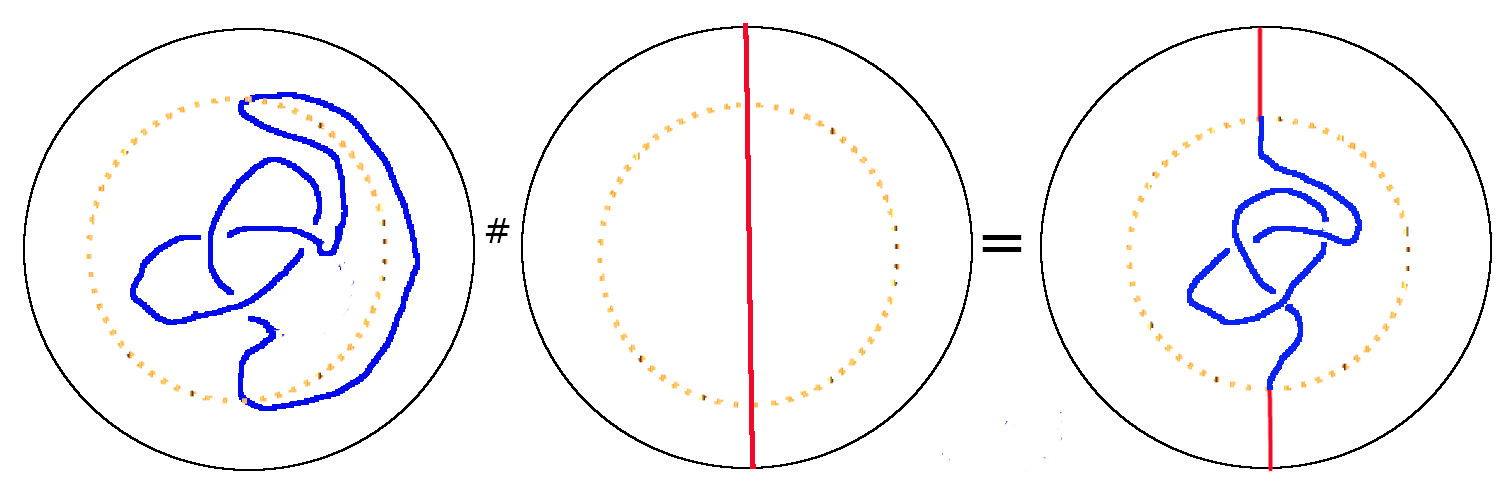}
\caption{Projectivization of an affine trefoil}
\end{figure}
\begin{definition}
A knot, $K$ in $\mathbb{R}P^3$ is called a \textbf{projectivized affine knot}, if there is a seperating sphere $S$ intersecting $K$ at exactly two points such that the corresponding $C\cap K$ is isotopic to a standardly embedded non-trivial cycle in $H_1(C,S)$.  
\end{definition}

\subsection{Tangles in the mapping cylinder}
Any knot in $\mathbb{R}P^3$ can intersect the mapping cylinders in a collection of curves. And if a knot, $K$, is not entirely contained in a mapping cylinder $C$, then $K\cap C$ is a tangle consisting of only arcs meeting the boundary sphere at exactly two points. They also can be linked within $C$. Refer to Figure 7. Thus the problem we are dealing with is the theory of tangles in the mapping cylinder. This tangle is in no way unique to the knot or to the separating sphere. In what follows we wish to study the tangles created by a knot in the mapping cylinders in $\mathbb{R}P^3$. If needed one can indeed push out a small arc from the knot outside the ball and add an extra component for the tangle. Hence it is in no way unique to the knot.\\
\begin{figure}
\begin{center}
\includegraphics{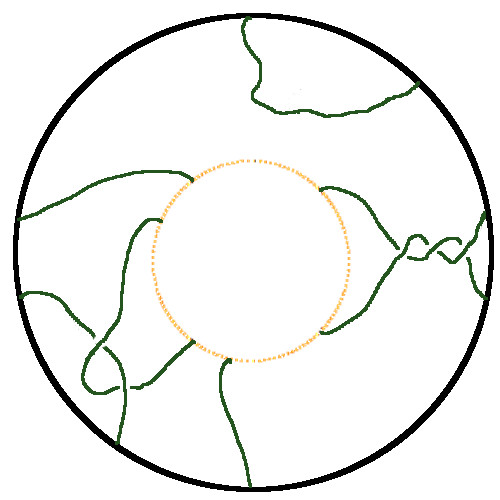}
\caption{A tangle in the mapping cylinder}
\end{center}
\end{figure}
Notice that $C$ is a 3-manifold with a boundary sphere and an asymptotic projective plane embedded in it. Whenever we say \say{the projective plane} in a mapping cylinder, it will always refer to the embedded asymptotic projective plane. And tangles, induced by knots, are composed of arcs with there boundary points on the sphere and interiors embedded in the interior of $C$. The arcs may or may not intersect the projective plane. The patterns of their intersections are definitely interesting to study. Such studies can give topological information on the knot. For example the $mod-2$ intersection number of the tangle and the projective plane is the homology class of the knot in the integral homology of $\mathbb{R}P^3$. Thus the $mod-2$ intersection number of the tangle induced by $K$ in $C$ with the projective plane is independent of the tangle and depends only on the knot.
Thus the projectivized affine knots, defined above, are knots with a simple tangle which consists of just one arc intersecting the projective plane at exactly one point.

\begin{figure}
\begin{center}
\includegraphics{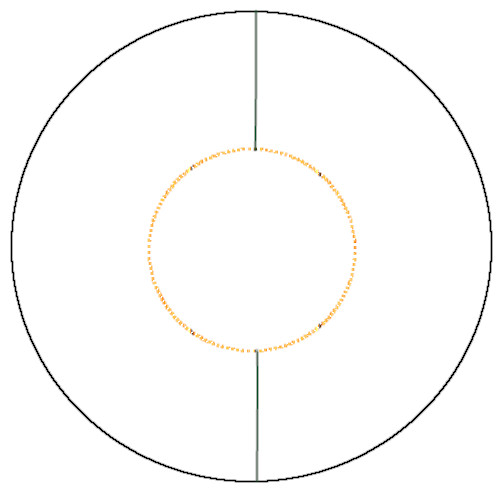}
\caption{Diagram of a class-1 arc in the mapping cylinder}
\end{center}
\end{figure}
 Since the affine knots are homologically trivial, one can see that the process of \say{projectivization} also changes a class-0 knot into a class-1 knot. And indeed this is a way of defining connected sums of a knot with the projective unknot which has a clear variant of connected sum with the affine unknot. But the procedure described above is limited to affine knots.
We would, here, like to define a general version of this process. And thereby obtain a classifying theme for all knots in $\mathbb{R}P^3$. 

\subsection{Residual tangles and space bending surgery}
Let $L^n$ represent the projective closure of $n$ disjoint non-parallel lines in $\mathbb{R}^3$. 
\begin{definition}
The tangle defined by $L^n$ in the complement mapping cylinder of some open ball centered at origin intersecting all lines in $L^n$ will be called the\textbf{ $n^{th}$ residual tangle} denoted by $T^n$.  
\end{definition}

\begin{figure}
\begin{center}
\includegraphics{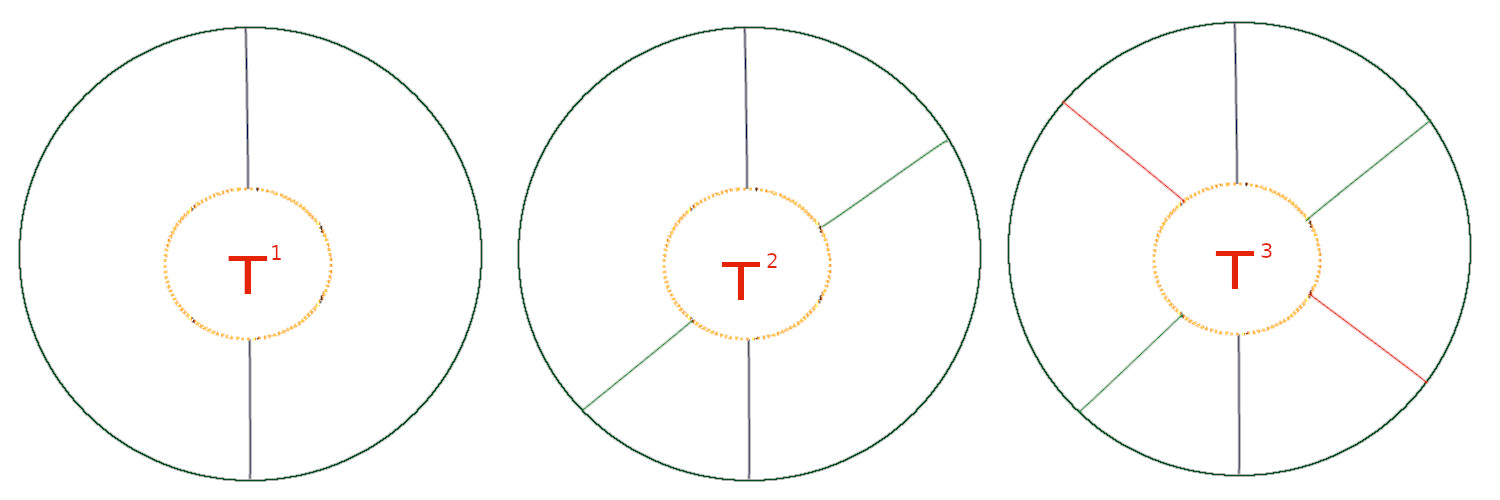}
\caption{Residual tangles}
\end{center}
\end{figure}
Notice that  $T^n$ is composed only of class-$1$ arcs in a mapping cylinder $C'$ intersecting the projective plane at exactly one point. Each arc has two boundary points on $S':=\partial C$. Hence the boundary of $T^n$ is a collection of $2n$ points on $S'$. It is easy to see that for every $n$ there is a unique such $T^n$ whose boundary has $2n$ points. \\

 Consider an affine knot $K$ in $\mathbb{R}P^3$ with a non-empty intersection with a disk $B$. Now the boundary of $B$ is a seperating sphere, $S$, and the complement of $B$ in $\mathbb{R}P^3$ is a mapping cylinder $C$. Notice that $C\cap K$ is a finite collection of arcs with boundary points on $S$. We may assume that they are all boundary parallel arcs, so that none of them are locally knotted. Thus $B\cap K$ is a tangle with some even number, say $2m$, boundary points on $S$.\\

Let $f: S\to S'$ be a diffeomorphism sending the boundary of $K\cap B$ to $\partial T^m$. 

\begin{definition}
Gluing $B$ and $C'$ along $f$ produces a copy of $\mathbb{R}P^3$,  with a new knot which is formed by gluing $K\cap B$ and $T^m$ along their boundary. We denote the knot thus obtained as $\Sigma(K,S,f)$. This procedure will be called as \textbf{\textit{space bending surgery.}}
\end{definition}

\begin{remark}  The philosophy behind the name is that, an affine knot is basically a classical knot in $\mathbb{R}^3$ with a new identity once we add a projcetive plane at infinity. Thus, through this surgery, we are connecting this knot with the projective plane at infinity, and we are making it \say{aware} of the change in the ambient space.
\end{remark}

Suppose $K$ is an affine knot and $S$ intersects $K$ at exactly two points. If, the arc in $C\cap K$ is also boundary parallel, then any $\Sigma(K,S,f)$ is a projectivization of $K$ as defined above. Thus we can clearly see that projectivization for a given affine knot is not unique. It depends on the choice of intersection of $K$ with $S$ and the diffeomorphism $f$.
The following theorem says, every knot in $\mathbb{R}P^3$ is isotopic to a knot of the form $\Sigma(K,S,f)$.

\begin{theorem}
{\textbf{\rm (The structure theorem):}} Any knot in $\mathbb{R}P^3$ can be obtained from an affine knot by a space bending surgery.
\end{theorem}

\textbf{Proof:} Let $J$ be a knot in $\mathbb{R}P^3$. Choose an asymptotic projective plane, $P$ such that $J$ intersects $P$ at $m$ points, transversely. let $C$ be a regular neighborhood for $P$. That is $C$ is a mapping cylinder considered to be an $I$-bundle over $P$. Closure of the complement of $C$ is a 3-ball $B$ with $S:= B\cap C= \partial B= \partial C$. \\
We may assume that, $J\cap C$ is just a collection of finitely many arcs all of whose boundaries are in $S$. We can also assume that all of these arcs are fibers of the $I$-bundle. All the other kind of arcs and local knots on them may be pushed inside $B$. Now, we can remove the arcs in $J\cap C$ with $m$ boundary parallel arcs in $C$. Then we can obtain an affine link, $K$.\\
Then it is obvious that $J$ is isotopic to $\Sigma(K,\partial B',f)$.\qed\\
The following diagram shows an example which demonstrates the dependence on the $f$ in the surgery.

\begin{figure}
\begin{center} 
\includegraphics{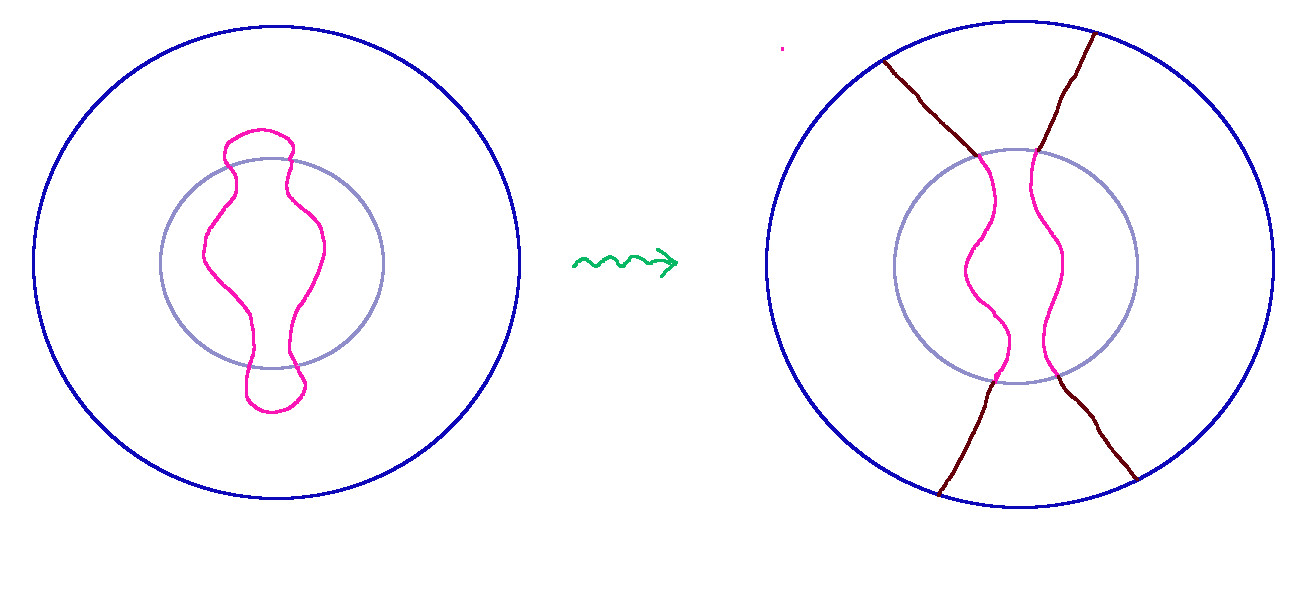}
\caption{Surgery on affine unknot producing affine unknot}
\end{center}
\end{figure}
\begin{figure}
\begin{center}
\includegraphics{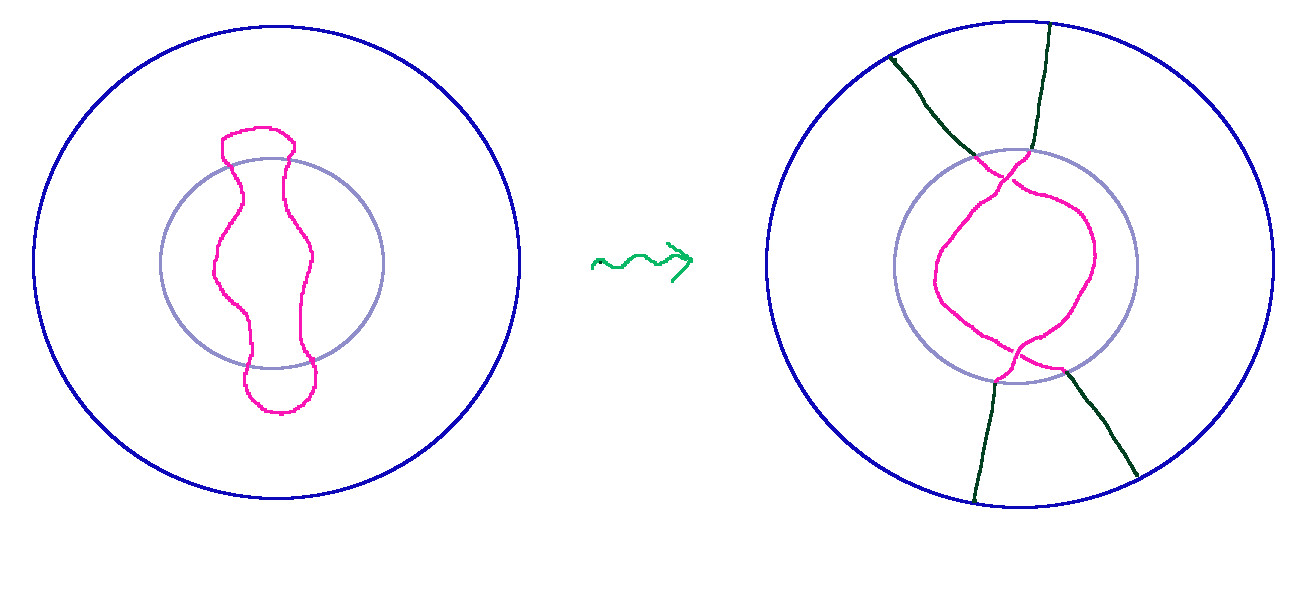}
\caption{Surgery on affine unknot producing a non-affine knot}
\end{center}
\end{figure}
The knots on the left are both obtained from the affine unknot by space bending surgery. The separating spheres are both the same. But by using two different $f$'s we can obtain two distinct knots. The one above is again an affine unknot and the one below is a non-affine class-0 knot.  \\
\begin{corollary}
For every knot in $\mathbb{R}P^3$, there exist a separating sphere such that in the corresponding separation, the tangle in the mapping cylinder is a residual tangle. 
\end{corollary}

\begin{definition}
Let $K$ be a knot in $\mathbb{R}P^3$. A separating sphere which induces a splitting such that the part of $K$ in the mapping cylinder is a residual tangle, is said to be a \say{residual sphere} for $K$.
\end{definition}

\subsection{Class changing surgery}
Space bending surgeries are reversible, and can be altered to produce many knots from a given knot. For example, let $S$ be a residual sphere for a knot $K$. Then in the induced splitting, $\mathbb{R}P^3:= B\cup C$, the tangle in $C$ is a residual tangle. Now by pushing an arc of $K$ from $B$ to $C$, such that the boundary points are diametrically opposite in $S$. Now remove this class-0 arc in $C$ and replace it by a class-1 arc in $C$. Thus we obtain a new knot (Figure 15). This can be repeated arbitrary number of times. This procedure is just a small variation to the space bending surgery. Each time it is performed, it changes the homology class of the knot. Thus we would refer to this procedure as \textbf{\say{class changing surgeries}}. The following picture demonstrate such a surgery on a non-affine knot.\\

\begin{figure}

\includegraphics{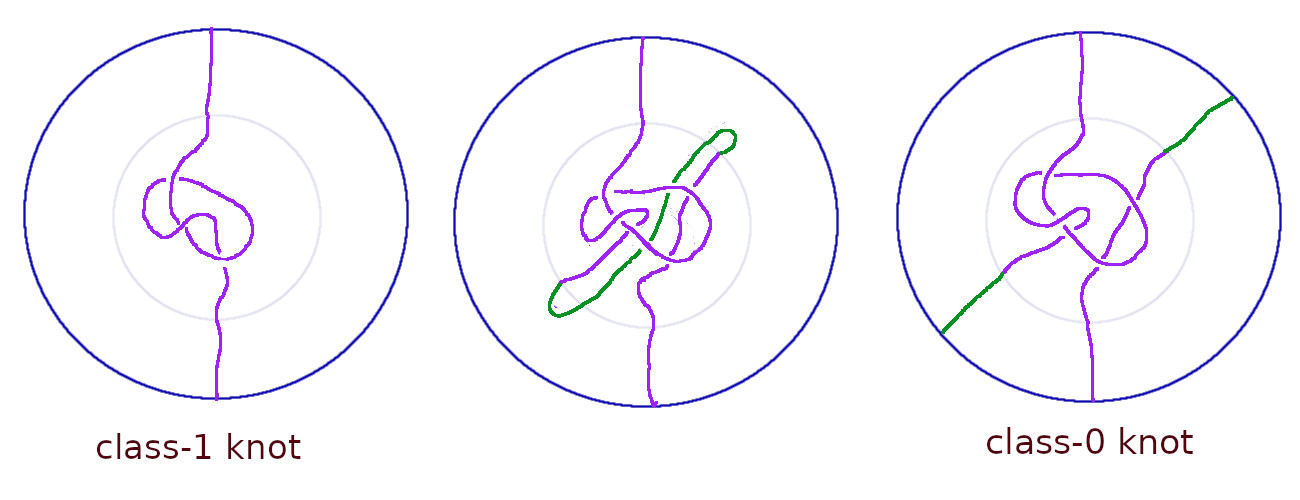}\\
\caption{Surgery on a projective trefoil producing an affine trefoil}

\end{figure}
Notice that the class-1 knot in the example above is a projcetivized affine trefoil. We may refer to it as a projective trefoil. It is easy to see that the structure theorem facilitates the surgery to be performed on arbitrary knots, rather than on just affine knots. This is because, when we are pushing out arcs of $K$ out of a ball $B$, if boundary of $B$ is a residual sphere, then it can be done in such a way that the arc of $K$ or the class-1 arc which is replacing it, would not be interlinked with the tangle outside.

\section{Genus for knots in $\mathbb{R}P^3$}
Let $S$ be a closed connected surface in $\mathbb{R}P^3$ containing a knot $K$. We know that $\widetilde{S}:=h^{-1}(S)$ is a closed surface in $S^3$ containing the link, $h^{-1}(K=)\widetilde{K}$. By removing any point in the complement of $\widetilde{S}$, we can obtain an embedding of $\widetilde{S}$ in $\mathbb{R}^3$. Every closed surface embedded in $\mathbb{R}^3$ is orientable, and hence, $\widetilde{S}$ is orientable. If $S$ is contained in some open 3-ball in $\mathbb{R}P^3$, that is, it is contained in the affine part of $\mathbb{R}P^3$, then $\widetilde{S}$ is disconnected. 
\begin{definition}
A closed surface $S$  containing $K$ in $\mathbb{R}P^3$ is said to be a \say{good surface for $K$} if the inverse image $\widetilde{S}:=h^{-1}(S)$ is connected.
\end{definition}
Notice that, if $S$ is a good surface for $K$, it may or may not be orientable. As an example, consider the affine unknot, $K_0$. Clearly, an asymptotic projective plane $P$ and a torus $T$ intersecting $P$ at only $K_0$, are both good surfaces for $K_0$. But $T$ is orientable, while $P$ is non-orientable. 

\begin{theorem}
Every knot in $\mathbb{R}P^3$ has a good surface. 
\end{theorem}
\textbf{Proof: } Let $K$ be a knot in $\mathbb{R}P^3$. Consider a diagram of $K$ on some aymptotic projective plane $P$. If $K$ is an unknot, (affine or projective) there is a diagram of $K$ with no crossing. The inverse image of $P$ is a great $2-sphere$ in $S^3$ and thus $P$ itself is a good surface for $K$. Suppose $K$ is knotted.  Then for each  crossing of $K$ we attach a handle on $P$ containing the crossing as shown in Figure 13. Thus,  we will obtain a surface, $S$ containing $K$. It is obvious that since $h^{-1}(P)$ is connected, $h^{-1}(S)$ is also connected. Thus $S$ is a good surface for $K$ and we are done. \qed \\
\begin{figure}
\includegraphics{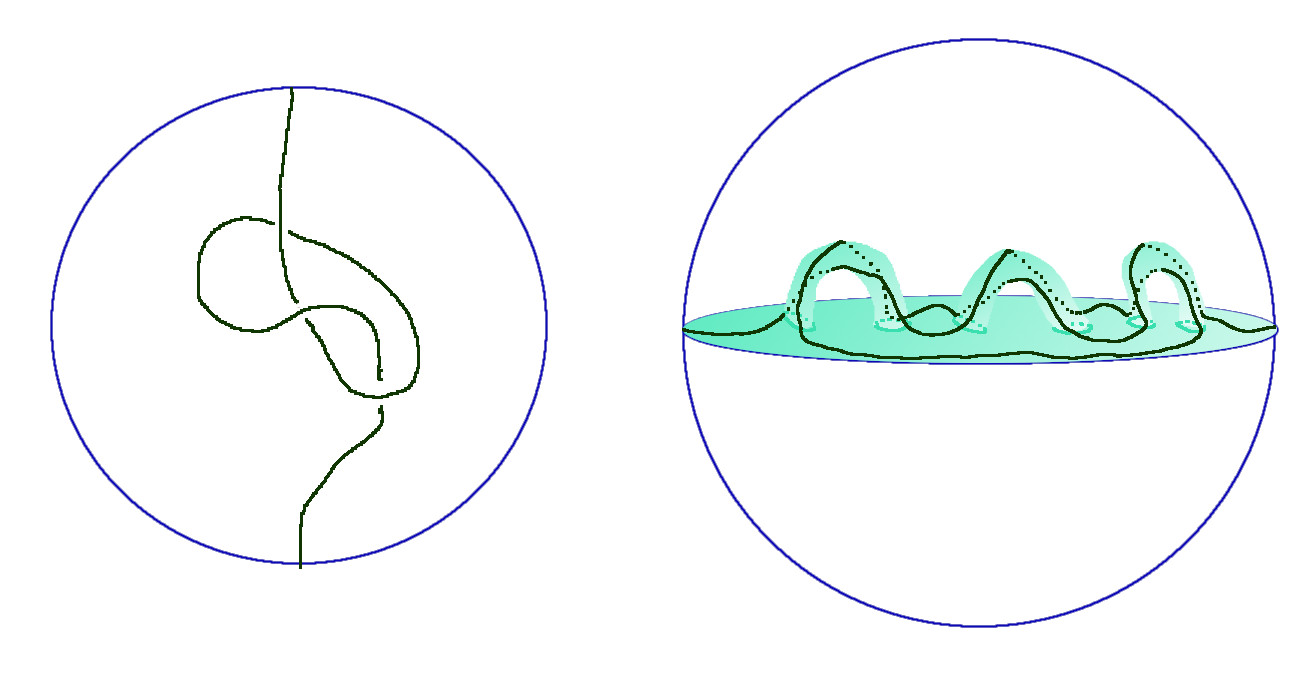}
\caption{A good surface containing projective trefoil}
\end{figure}

\begin{remark} The good surfaces obtained by the above construction are all non-orientable.  However we will see later that there are orientable good surfaces for some knots $K$. The authors are not aware when exactly a given knot has an orientable good surface. 
\end{remark}

\begin{definition}
The minimum genus of an inverse image surface varied over all the good surfaces for $K$ will be called the \say{genus of $K$}. 
\end{definition}

\begin{theorem}
A knot in $\mathbb{R}P^3$ is an unknot if and only if it has genus 0. 
\end{theorem}
\textbf{Proof:} By definition, both unknots are embedded in an asymptotic projective plane. Since any asymptotic projective plane is covered by a 2-sphere in $S^3$, both the unknots have genus 0. \\
Conversely, if a knot has genus 0, then it has a good surface which is covered by a 2-sphere in $S^3$. By definition, this sphere has to be isotopic to a great 2-sphere in $S^3$. Thus the good surface is an asymptotic projective plane. And hence the knot is an unknot. 
Thus we are done. \qed

As a simple consequence of the theorem above, we derive the following \say{non-cancellation property} of the space bending surgery.

\begin{corollary}
A space bending surgery on a non-trivial link, will always yield a non-trivial link. 
\end{corollary} 
\textbf{Proof:} Let $K'$ be a link obtained by performing a space bending surgery on a non-trivial link $K$. Let $B$ and $C$ be the ball and mapping cylinder used in the surgery. If $K'$ is an unlink, then by theorem 4.3, it follows that it has a good surface $\Sigma$ which is an asymptotic projective plane. Let $F:\mathbb{R}P^3\times[0,1] \to \mathbb{R}P^3$ be an ambient isotopy such that $F_0=Id$ and $F_1(K')\subset \Sigma$. This implies that, $B\cap K=B\cap K'$ is a tangle embedded in the region, $X:=F_1(B)\cap \Sigma$. We may assume that this $F_1(B)$ intersects $\Sigma$ transversally. $F_1(B)$ is homeomorphic to a 3-ball and hence $X$ should be a finite collection of 2-disks in $\Sigma$. Thus, all arcs in $X$ are unknotted. Which implies that all arcs in the tangle $B\cap K$ also were unknotted. We know that $K\cap C$ is a collection of boundary parallel arcs which are also unknotted. Hence, $K$ should be an unlink which contradicts the initial hypothesis. Thus, $K'$ cannot be an unlink and we are done. \qed

\begin{theorem}
All non-trivial class-$1$ knots have genus 1.  
\end{theorem}
\textbf{Proof: } Let $K$ be a class-$1 $ knot in $\mathbb{R}P^3$. Let $U$ be a tubular neighborhood of $K$. The boundary of $\bar{U}$ is a torus. By choosing an arbitrary non-vanishing vector field, $X$ on $K$ and pushing $K$ into $T:=\partial \bar{U}$ along $X$ we can obtain a parallel of $K$ isotopic to $K$ inside $T$. Now the inverse image of $U$ is a tubualar neighborhood of $\widetilde{K}$ and hence is a solid torus in $S^3$. The boundary of this tubular neighborhood is a torus containing the inverse image of the parallel of $K$ contained in $T$. Thus $T$ is a good surface for $K$ (or its parallel). And since $K$ is non-trivial, the genus of $K$ is $1.$ Hence we are done. \qed

Now we will study the properties of genus for affine knots. Notice that any affine knot, has two lifts in $S^3$ which are unlinked. Any of these lifts, being a knot in $S^3$ has a genus. We would like to compare, for an affine knot, the genus of a lift in $S^3$ and the genus the knot in $\mathbb{R}P^3$. \\
Consider an orientable surface, $S$, contained in an affine region of $\mathbb{R}P^3$. Let $P$ be an asymptotic projective plane disjoint from it. By removing a small disk from both $S$ and $P$, and connecting them by a tube, we can always obtain a surface in $\mathbb{R}P^3$ which is homeomorphic to the connected sum, $P\#S$. The inverse image of this surface would be a great sphere (the inverse image of $P$) attached with two copies of $S$ with tubes, one \say{inside} and one \say{outside}. This is a connected surface. Thus for any knot $K$ contained in a surface $S$ in an affine region, the connected sum of $S$ with a projective plane is always a good surface for $K$.  

\begin{lemma}
Every non-trivial affine knot has genus strictly greater than $1$. 
\end{lemma}
\textbf{Proof: }Let $K$ be a non-trivial affine knot. If $K$ had genus 0, then it is an unknot by Theorem 4.3, hence cannot be a non-trivial knot. Thus $K$ cannot be of genus 0. Now suppose $K$ had genus $1$. Then it has a good surface $G$, whose inverse image is a torus in $S^3$, which is a two sheeted cover for $G$. Then, $G$ is Klein's bottle or a torus. Since the Klein's bottle cannot be embedded in $\mathbb{R}P^3$, see \cite{bredon}, $G$ is a torus.  \\

We know that, since $K$ is an affine knot, there exist an asymptotic projective plane, $P$, disjoint from it. If $P\cap G=\emptyset$ then, $h^{-1}(G)$ is disconnected which contradicts the fact that $G$ is a good surface. Hence $G$ intersects $P$ non-trivially. We may assume the intersection of $P$ with $G$ is transversal and thus is a collection of circles. Suppose one of these circles, $\sigma$ is contractible in $G$. Then $\sigma$ is homologically trivial in $G$, $P$ and $\mathbb{R}P^3$. Then it bounds a 2-disk $D\subset G$ and a 2-disk $D_1\subset P$. Notice that, $D_1\cap G$ is a subcollection of the collection of circles $P\cap G$. This subcollection contains $\sigma$. If any other circles belong to this subcollection, they are embedded in the interior of $D_1$. Choose a regular neighborhood of $\sigma$ in $P$, say $U$. Then $U$ is homeomorphic to a cylinder and $\partial \overline{U}= \sigma_1\cup \sigma_2$, two disjoint circles. Exatly one of these circles are disjoint from $D_1$, say $\sigma _2$. Now choose a regular neighborhood $W$ of $D$ in $\mathbb{R}P^3$. $\overline{W}$ is a trivial $I$-bundle over $D$. Choose a homeomorphism, $f:D\times [0,1]\to \overline{W}$. We may assume that $U$ and $W$ are chosen in such a way that $U=f(\sigma\times (0,1))$. Then $\partial \overline{W}= U\cup f(\{0,1\}\times D)$. Let $D_2$ and $D_2'$ be the two disks in $f(\{0,1\}\times D)$. Exactly one of the disks, say $D_2$ is bounded by $\sigma_2$. Now, we may remove the disk bounded by $\sigma_2$ inside $P$ and replace it with $D_2$. Clearly this gives a new projective plane which is isotopic to $P$ and it doesn't intersect $G$ at any of the circles in $D_1\cap G$. We may call this projective plane also as $P$ and think of it as a modification of the old $P$ with a contractible circle in its intersection with $G$ \say{removed}. \\
 
  We can use this procedure to remove all the contractible circles in $P\cap G$. Thus we may assume all circles in $P\cap G$ are essential in $G$. Hence $G\setminus P$ is a collection of cylinders. Since $K$ is a knot, it is contained in exactly one cylinder, say $C$, in $G\setminus P$. Being disjoint from $P$, $C$ is contained in an affine region of $\mathbb{R}P^3$, the complement of $P$. Let $V$ be a regular neighborhood of $P$ disjoint from $K$. Then $C':=C\setminus V$ is a closed cylinder with two boundary circles. Notice that $K\subset C'$. Since $K$ has no self intersections, $K$ with some orientation represents a generator for $H_1(C',\mathbb{Z})$. But we know that a boundary circle of $C'$, say $\alpha$, is another representative for a generator of $H_1(C',\mathbb{Z})$. Since $\alpha$ is contained in the complement of $P$ it is an affine unknot. Then, $K$ is homologous to $\alpha$ in $C$ and together they bound a cylinder inside $C$. By shrinking this cylinder we can construct an isotopy of $K$ to $\alpha$. Thus $K$ is isotopic to the affine unknot. Which is a contradiction. Thus the genus of $K$ has to be strictly greater than $1$. \qed\\

\begin{theorem}
The genus of any non-trivial affine knot in $\mathbb{R}P^3$ is $2. $
\end{theorem}
\textbf{Proof: } Let $K$ be an affine knot,  contained in the affine part of $\mathbb{R}P^3$. Consider a torus (possibly knotted) containing it in the affine part and take a connected sum with a projective plane. See Figure 14. Thus we can get a good surface for $K$ in $\mathbb{R}P^3$. Then the genus of the inverse image of this good surface is $2$. This implies that the genus of any affine knot, is less than or equal to $2$. Then by Lemma 4.5, the genus of $K$ is exactly $2.$ \qed \\

\begin{figure}
\begin{center} 
\includegraphics{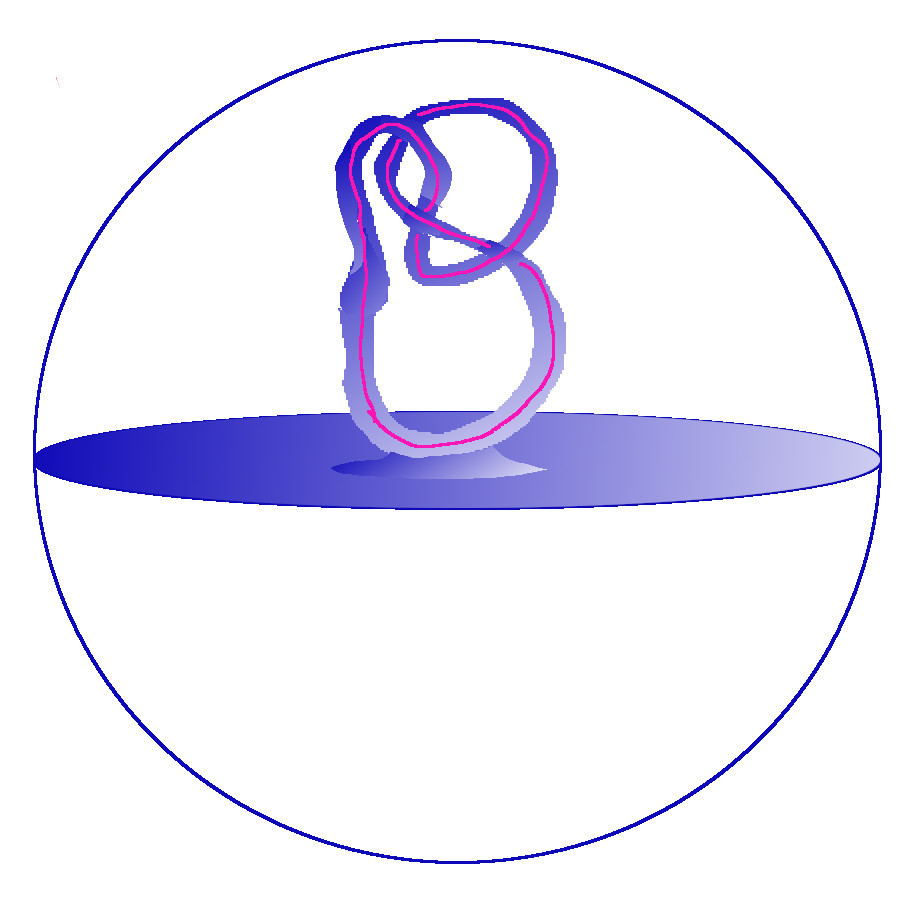}
\caption{A minimal good surface for affine trefoil}
\end{center}
\end{figure}

\begin{corollary}
All knots of genus $1$ are non-affine.
\end{corollary}
Just like knots, any subset of $\mathbb{R}P^3$ has a notion of being affine or non-affine.\\
\begin{definition}
We will call a torus $T^2$ in $\mathbb{R}P^3$ to be class-$1$ if it is the boundary of a tubular neighborhood of a class-$1$ knot. In particular the boundary of the tubular neghbourhood of a projective unknot will be called a standard class-1 torus. 
\end{definition}
Clealy the inverse image of such a torus under $h$ is a connected torus in $S^3$. Thus for a knot embedded in a class-$1$ torus, it is a good surface. Since every knot of genus $0$ is an unknot, a non-trivial knot embedded on a class-$1 $ torus has genus $1.$\\
\begin{figure}
\begin{center}
\includegraphics[scale=1]{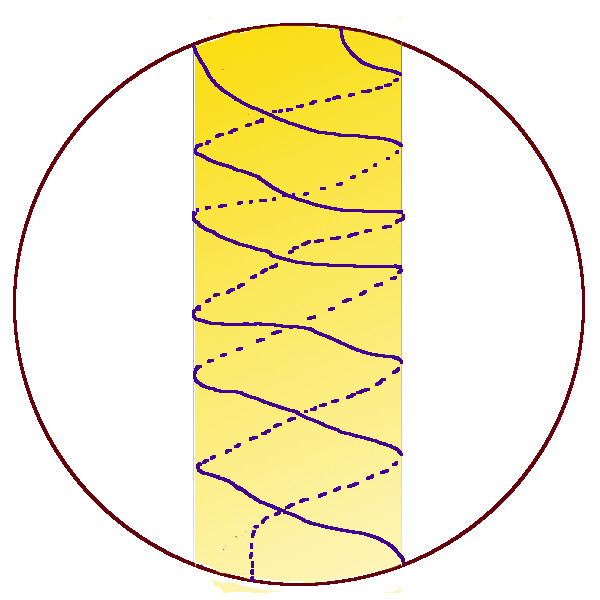}
\caption{A class-0 non-affine knot of genus 1 on a class-1 torus}
\end{center}
\end{figure}
The first homology of the torus, $H_1(T^2,\mathbb{Z})\approx \mathbb{Z}\oplus \mathbb{Z}$. Any knot on the torus $T^2$ represents a homology class of the form $(p,q)$ where $p$ and $q$ are coprime. Choose a ball model of $\mathbb{R}P^3$ as quotient of a 3-ball $D^3$. Consider knot shown in Figure 15 as a diagram drawn on a 2-disk $D^2$. Let $d$ be a diameter of $D^2$, then $\pi(d)$ is a projective unknot. Choose a tubular neighborhood $U$ of $\pi(d)$. Without loss of generality, we may assume, the image of $U$ under $\delta$ is a region in $D^2$ which bounded by two parallel chords of $\partial D^2$. See the shaded region in Figure 15. Clearly $T:=\partial \overline{U}$ is a standard class-1 torus intersecting the asymptotic projective plane $P:=\pi(\partial D^3)$ at exactly one meridinal circle, say $m$, of the solid torus $\overline{U}$. Let $K$ be a $(p,q)$ knot in $T$ which is the quotient of a torus braid embedded in the cylinder $\pi^{-1}(T)$. Figure 15 gives an example of a braid whose quotient is the $(2,5)$ knot in $T$. All points in $K\cap P$ are contained in $m$. We can choose $K$ in such a way that there are exactly $p$ points in $K\cap P$. Thus $K$ is class-0 whenever $p$ is even. Hence if it is non-trivial, then it must be non-affine by corollary 4.8. \textbf{Therefore, there are infinitely many class-$0$ non-affine knots with genus $1$.}
\begin{figure}
\begin{center}
\includegraphics{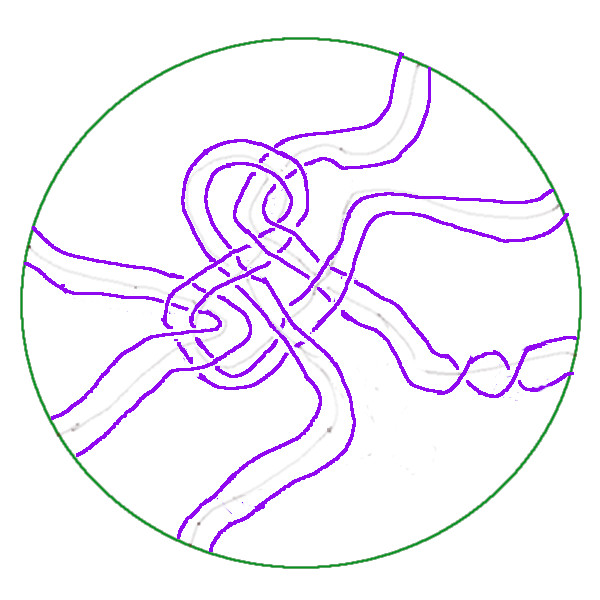}
\caption{A class-$0$ non-affine knot on a knotted class-$1$ torus}
\end{center}
\end{figure}

The knot shown in Figure 16, is non-affine, since the inverse image in $S^3$ is a two component link with linking number 2 between its components. Clearly this knot, is boundary of a band in $\mathbb{R}P^3$. This band can be thickened to be a class-$1$ torus on which the knot represents a $(2,5)$ cycle. Thus it has genus $1.$ The core of the solid torus bounded by this torus is class-$1$ knot intersecting an asymptotic projective plane at infinity at 3 points. 
\begin{conjecture}
There exists class-$0$ non-affine knots, with genus not equal to $1. $
\end{conjecture}
As a support for the conjecture, in Figure 20, we present a diagram of a knot, which we suspect to be of genus $3. $ The surface on which the knot is lying is good and its inverse image is a genus-$3 $ surface in $S^3$.
\begin{figure}
\begin{center}
\includegraphics{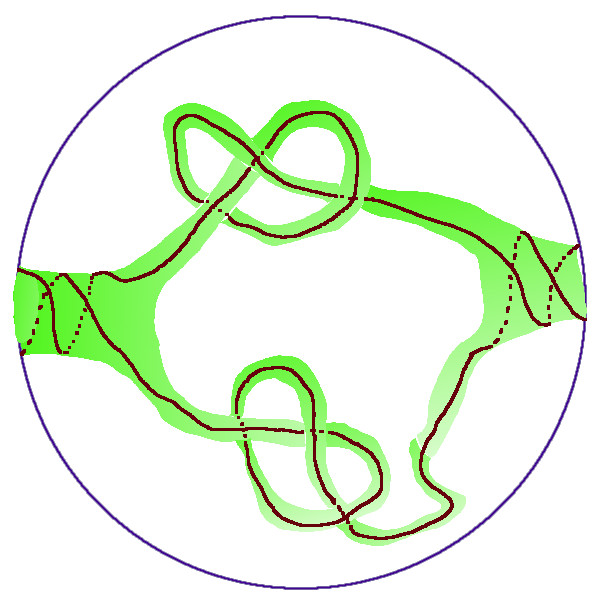}
\caption{A non-affine knot on a good surface}
\end{center}
\end{figure}
\section{\textbf{Companionship and Projective inessentiality}}
The following is a definition of companionship for knots in $\mathbb{R}P^3$ directly generalized from that of $S^3$. A subset of a solid torus is said to be geometrically essential, if every meridinal disk intersects it. 

\begin{definition}
A knot $K_1$ in $\mathbb{R}P^3$ is called a companion for a knot $K_2$ if there exists a tubular neighborhood, $V$ of $K_1$ containing $K_2$ such that $K_2$ is geometrically essential in $V$.
\end{definition}

Notice that $\mathbb{R}P^3$ is diffeomorphic with $L(2,1)$, a lens space. A tubular neighborhood of a projective unknot, is a solid torus  whose complement is also a solid torus. Thus the boundary of a tubular neighborhood of a projective unknot induces a Heegaard splitting for $\mathbb{R}P^3$. The meridian of each of the solid tori appearing in the splitting is glued to the $(2,1)$ curve on the boundary of the other one in order to construct $L(2,1)$. The characteristic curve on the Heegaard splitting surface, bounds a 2-disk in one solid torus and a Mobius strip on the other. They together glue to form an asymptotic projective plane in $\mathbb{R}P^3$. 

\begin{theorem}
 The projective unknot is a companion of every knot $K$, in $\mathbb{R}P^3$. 
\end{theorem}
\textbf{Proof:} Choose a ball model of $\mathbb{R}P^3$ presented as quotient of a ball $D^3$. Consider a diagram of $K$ on a disk $D^2$, as in Section 2. Choose a separating circle $\sigma$ which separates the diagram into the diagram of a residual tangle in the annulus $A$ and the diagram of the complementary tangle in the disk $D\setminus A$.\\
\textbf{Case 1:} Suppose $K$ is a non-affine knot. Then it has a non-empty residual tangle in $A$. Now, as shown in Figure 18, choose a pair of half disks intersecting the boundary circle of $D^2$.
\begin{figure}
\begin{center}
\includegraphics{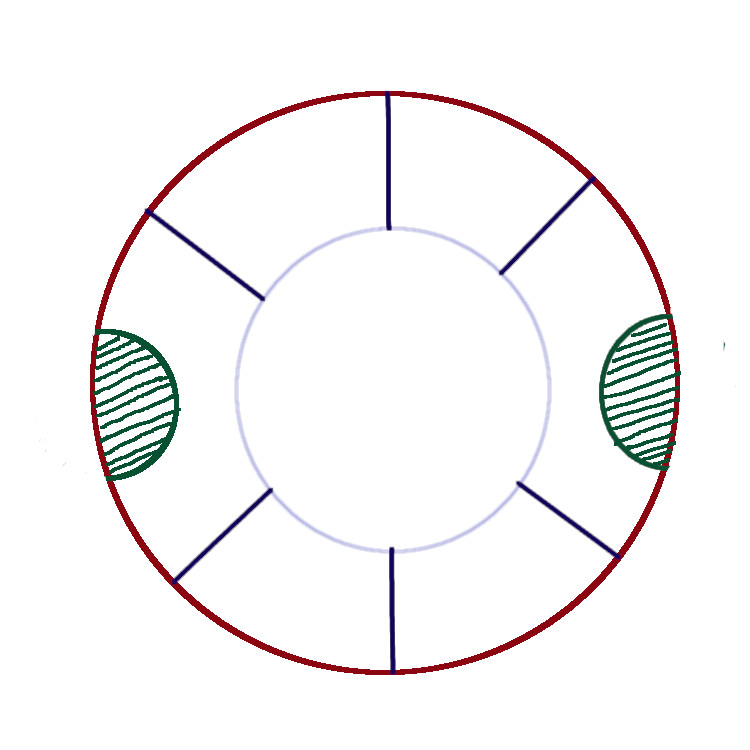}
\caption{The shaded region is a 2-disk in the quotient $\mathbb{R}P^2$}
\end{center}
\end{figure}
The two half disks represents a $2$-disk, say $B$ in the quotient projective plane,$P:=\pi(D^2)$. Let $\overline{N}\in \mathbb{R}P^3$ be the point on which $\delta$ is not defined, as in Section 2. Then $V:=\pi^{-1}(B)\cup\overline{N}$, is diffeomorphic to a solid torus pinched at a point. Refer to Figure 19.
\begin{figure}
\begin{center}
\includegraphics{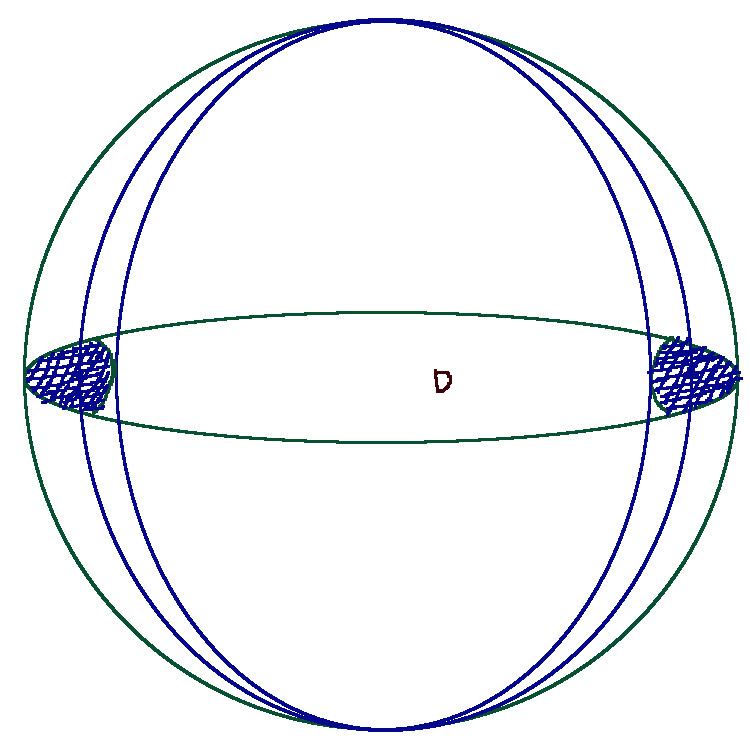}
\caption{The pull back of a disk under the projection  representing a solid torus pinched at a point}
\end{center}
\end{figure}
 We may choose a small enough ball neighborhood, say $U$ of $\overline{N}$ such that $U\cup V$ is homeomorphic to a solid torus. It can be smoothed out to be a solid torus, $V_2$, in $\mathbb{R}P^3$ disjoint from $K$. Notice that, if $x$ is the center of the disk $B$, then $\delta^{-1}(x)\cup \overline{N}$ is a projective unknot embedded in the interior of $V_2$. Thus $V_2$ induces a Heegaard splitting of $\mathbb{R}P^3$. The complement, $V_1:=\mathbb{R}P^3\setminus V_2$, is an open solid torus containing $K$.\\
 
  Then, $\mathbb{R}P^3= V_1\bigcup V_2$ is a Dehn surgery description of $\mathbb{R}P^3$ as obtained from $S^3$ by a $(2,1)$ surgery on the unknot. Notice that the core circle of $V_1$, say $L$ is a projective unknot in $\mathbb{R}P^3$ and $V_1$ is a tubular neighborhood of $L$. Now since two 1-dimensional submanifolds in $V_1$ cannot intersect transversally, we may assume that $L$ is disjoint from $K$. Now any meridinal disk in $V_1$ shares its boundary with a Mobius strip in $V_2$ and together they form an asymptotic projective plane. Thus if $K$ was disjoint from a meridinal disk then it affine, which is a contradiction to the initial assumption. Thus $L$ is a companion of $K$. \\
  
\textbf{Case 2:} Suppose $K$ is an affine knot. Consider a diagram of $K$ on a 2-disk $D^2\subset D^3$ with a separating circle $\sigma$ as in the previous case. Let $A$ and $B$ be the corresponding open annulus and 2-disk in splitting induced by $\sigma$, that is, $D^2\setminus \sigma=A\sqcup B$. By definition, we may assume the diagram is contained entirely inside $B$. There is a unique 3-ball $B^3\subset D^3$ such that $B^3\cap D^2= B$. The quotient map $\pi$ is injective on $B^3$ and $\pi(B^3)\approx B^3$. We will refer to $\pi(B^3)$ also as $B^3$. Then $U:=D^3\setminus \overline{B^3}\approx S^2\times (0,1)$ is the unique $(0,1)$-bundle over $S^2$ such that, $U\cap D^2=A$. Notice that $\delta(B^3)=B$ and hence $K\subset B^3$. Let $C:=\pi(U)$ which is a mapping cylinder containing the asymptotic projective plane $P:=\pi(\partial D^3)$. Now by pulling out a boundary parallel arc of $K$ out of $B^3$ into $C$ and moving it through $P$ back into $B^3$ we may create a non-empty residual tangle, $T$ for $K$ with two arcs in $C$. Then the new diagram on $D^2$, $\delta(T)$ is a collection of four arcs in $A$ and $A\setminus\delta(T)$ is four disjoint regions. By isotopically moving $K$ if required, we may assume that there is diameter $\gamma$ of $D^2$ disjoint from the diagram of $K$. Then $\gamma$ intersects exactly two of these regions which are diametrically opposite. Choose an open 2-disk $E$ in the Mobius strip $M:=\pi(A)$ such that the two half disks in $\pi^{-1}(E)$ are contained in the two regions in $A\setminus\delta(T)$ which are disjoint from $\gamma$. Now choose an open solid torus $V$ containing $\delta^{-1}(E)$ as in previous case. Then $\overline{V}$ is clearly disjoint from $K$. Let $x$ be the central point of $E$. The core of $V$, $\delta^{-1}(x)\cup \{\overline{N}\}$ is a projective unknot. Then $V':=\mathbb{R}P^3\setminus \overline{V}$ is a solid torus and $\mathbb{R}P^3=\overline{V}\cup\overline{V'}$ is a Heegaard splitting. $V'$ contains $K$ and the projective unknot $\pi(\gamma)$. Without loss of generality, we can assume that $V'$ is a tubular neighborhood for $\pi(\gamma)$.\\

Let $F$ be a meridinal disk of $V'$ and let $\alpha:=\delta(F)$ a proper arc with its boundary points in $\partial\overline{E}$ intersecting $\pi(\gamma)$ transversally at one point. Suppose $K$ is disjoint from $F$. Notice that $F$ is contained in an asymptotic projective plane $Q$ disjoint from $K$. Let $P':=\delta^{-1}(d)\cup \{\overline{N}\}$, another asymptotic projective plane disjoint from $K$. We may assume that $Q$ and $P'$ intersect transversally at a projective line. . Then, $X:=\mathbb{R}P^3\setminus (P'\cup Q)$ is a disconnected space with exactly two components. Also notice that, $Y:=V'\cap X$  and $\delta(Y)=\pi(D^2)\setminus (E\cup\pi(\gamma)\cup\alpha)$ both have two components. Since the diagram of $K$ intersects both components of $\delta(Y)$ non-trivially $K$ should intersect both components of $Y$ non-trivially. But since $K$ is connected and is disjoint from both $P'$ and $Q$, it should lie in exactly one connected component of $X$. Which is a contradiction. Thus, $K$ interesects every meridinal disk of $V'$. Thus $K$ is geometrically essential in $V'$ and $\pi(\gamma)$ is a companion for $K$. Hence we are done. \qed 

\begin{figure}
\begin{center}
\includegraphics[scale=0.7]{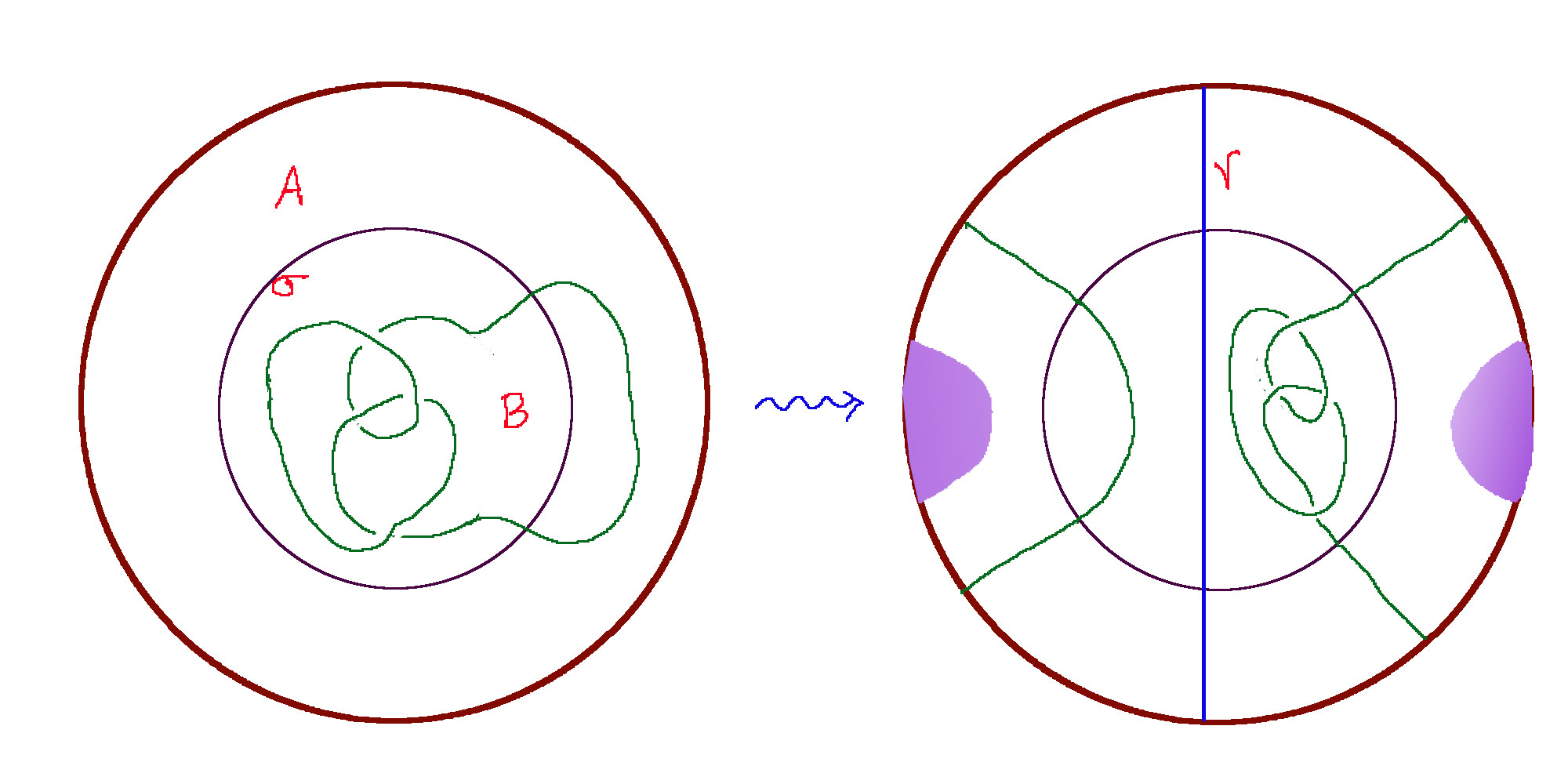}
\caption{The shaded region is $\pi^{-1}(E)$ }
\end{center}
\end{figure}

\begin{definition}
A Mobius strip, $M$, neatly embedded in a solid torus $V$, such that $\partial M$ represents the $(2,1)$ class in $H_1(\partial V,\mathbb{Z})$ will be called an \textbf{untwisted Mobius strip}. 
\end{definition}

\begin{definition}
Let $K$ be a knot in $\mathbb{R}P^3$. $K$ will be called \say{\textbf{projectively inessential}} if it is contained in a solid torus $V$ from a Heegaard splitting such that there exists an untwisted Mobius strip in $V$ disjoint to $K$.
\end{definition}

\begin{theorem}
A knot in $\mathbb{R}P^3$ is affine if and only if it is projectively inessential.
\end{theorem}
\textbf{Proof: (If:)}Let $K$ be a projectively inessential knot. By definition, there exist a solid torus $V$ containing $K$ whose complement $U:=\mathbb{R}P^3\setminus \interior{V}$ is also solid torus. $\mathbb{R}P^3=V\cup U$ is a Heegaard splitting. And there is an untwisted Mobius strip $M$ in $V$ disjoint from $K$. The circle $\partial M$ is a characteristic curve for the splitting and hence is a meridian for $U$ which bounds a disk $D$ in $U$. Then $D\cup M$ is an asymptotic projective plane and it is disjoint from $K$. Thus $K$ is an affine knot.   \\
\textbf{(Only if:)} Suppose $K$ is an affine knot in $\mathbb{R}P^3$. Choose an asymptotic projective plane $P$ disjoint from $K$. There exists a ball model for $\mathbb{R}P^3$ where $P$ is the image of an equitorial disk, $D$. Choose a disk $D'$ perpendicular to $D$ and consider a diagram of $K$ on $D'$ given by a projection with $N$ and $N'$ on $\partial D\setminus D'$. The image of $P$ under the projection is a projective line $L$, disjoint from the image of $K$. Choose a residual sphere for $K$ which will intersect $P'$ on a circle. Choose a 2-disk on $P'$ disjoint from the diagram of $K$ and $L$. Then the inverse image of this disk can be modified in to a solid torus $V$, whose complement solid torus $U$ contains both $K$ and $L$. Since $K$ is tame, $V$ can be chosen in a way that, it is disjoint from $K$ and $\partial V$ is a transversal to $P$. Clearly $V\cap P$ is a 2-disk, hence $U\cap P$ is an untwisted Mobius strip. Thus $K$ is projectively inessential. And we are done. \qed\\

\par This theorem provides a purely geometric criterion for a knot to be affine. We also wish here to mention a theorem announced by J.Viro and O. Viro in \cite{viro}

\begin{theorem}[J.Viro-O.Viro]
A knot $K$ is affine if and only if $\pi_1(\mathbb{R}P^3\setminus K)$ has an element of order 2. 
\end{theorem}
This is a purely algebraic criterion for a knot to be affine. Thus we have many techniques to detect whether a given knot is affine.

\section{acknowledgement}
The authors are thankful to the referee for valuable comments.

\end{document}